\documentclass[a4paper,10pt]{amsart}

\usepackage{amstext}
\usepackage{amsmath}
\usepackage{ifthen}
\usepackage{amscd}
\usepackage{amssymb}
\usepackage[english]{babel}
\usepackage[T1]{fontenc}
\usepackage[utf8]{inputenc}
\usepackage[all]{xy}
\usepackage{graphicx}
\usepackage{enumerate}
\usepackage{xspace}
\usepackage{pdfsync}
\usepackage{epic}

\newcommand{\wtx}{\tilde{X}}

\newcommand{\xu}{\widetilde{X^{un}}}

\newtheorem{prop}{Proposition}[section]
\newtheorem{coro}[prop]{Corollary}
\newtheorem{lem}[prop]{Lemma}
\newtheorem{rem}[prop]{Remark}
\newtheorem{exe}[prop]{Example}
\newtheorem{defi}[prop]{Definition}
\newtheorem{theo}[prop]{Theorem}
\newtheorem{conj}{Conjecture}

\newtheorem{ques}[prop]{Question}
\newtheorem{theoi}{Theorem}

\newcommand{\Q}{\mathbb Q}
\newcommand{\C}{\mathbb C}
\newcommand{\N}{\mathbb N}

\newcommand{\Z}{\mathbb Z}
\newcommand{\X}{\mathcal X}

\newcommand{\wt}{\widetilde}

\newenvironment{prv}{Proof:}{$\Box$}

\begin{document}

\title[Examples of Orbifold K\"ahler Groups]{Orbifold K\"ahler Groups and the Shafarevich Conjecture for Hirzebruch's covering surfaces with equal weights}

\date{\today}
\author{Philippe {Eyssidieux}}
\dedicatory{To Ngaiming Mok, for his 60th birthday. }

\begin{abstract}
This article is devoted to examples of (orbifold) K\"ahler groups from the perspective of the so-called Shafarevich conjecture on holomorphic convexity.
It aims at pointing out that every quasi-projective complex manifold with an \lq interesting\rq \ 
fundamental group gives rise to  interesting instances of this long-standing open question.

Complements of line arrangements are one of the better known classes of quasi-projective complex surfaces with an interesting 
fundamental group.  
We solve  the corresponding instance of the Shafarevich conjecture partially giving a proof that the universal covering surface of a
Hirzebruch's covering surface with equal weights is holomorphically convex. 

The final section reduces the Shafarevich conjecture to a question related to the Serre problem.
\end{abstract}

\maketitle

This article is the result of my effort  to 
understand the candidate counterexamples  given in \cite{BK} to the Shafarevich conjecture that the universal covering space of a complex projective manifold 
should be holomorphically convex. 
Being unable to settle them,  I
 recasted them in stack-theoretic terms and generalized the construction
to generate a wealth of similar examples where the Shafarevich conjecture leads to pretty
non-trivial group theoretic statements, some of which are in principle much simpler than the conjecture in \cite{BK}. 

In view of \cite{EKPR}, a counterexample to the Shafarevich conjecture should have a pretty complicated fundamental group since it
is not a linear group, i.e.: it could not be realized as a subgroup
of a matrix group. Furthermore, when one knows that the fundamental group has finite dimensional complex linear representations with infinite image,
the problem becomes tractable, although non trivial. 

On the other hand, if there are non-linear (actually non-residually finite) complex-projective groups ($\pi_1$  of complex projective manifold) \cite{Tol}, 
it is not known whether there exists a (possibly residually finite) infinite complex-projective group with no
finite dimensional complex linear representations with infinite image. 
Hence non trivial instances of the Shafarevich conjecture  are to be constructed from non trivial examples of complex-projective (or K\"ahler) groups. 

Now, it is much easier to construct {\em orbifold} K\"ahler groups than K\"ahler groups themselves.
They appear as quotients of fundamental groups of quasi-K\"ahler manifolds
and the first section will begin by reviewing this classical construction with more focus on its stack-theoretic aspects than in the litterature.
The outcome of this construction is that when one has an interesting fundamental group of a quasi-projective complex manifold, 
then one has a wealth of potentially interesting orbifold K\"ahler groups\footnote{
If residually finite, an orbifold K\"ahler group is K\"ahler \cite{vjm} and it is not clear whether the two classes are the same. 
This is actually one of the main difficulties in proving that the classes of K\"ahler and complex-projective 
groups coincide.}. Propositions \ref{nonunif} and \ref{develo} are the main original results in this section and may be 
stated as a Theorem answering questions asked in \cite{vjm}: 

\begin{theoi}
Every orbifold K\"ahler group is the fundamental group of a developpable compact K\"ahler orbifold.
There exists a non uniformizable developpable compact K\"ahler orbifold.
\end{theoi}

A word is due on the stack-theoretic perspective adopted here. Among its numerous advantages over 
Thurston's approach to orbifolds  advocated by \cite{art:lerman2010}, I want to emphasize the clear definition of a map of 
stacks, of substacks of a given stack, more generally of (2-)fibered products of stacks and add that there is now enough fundational material 
on the algebraic topology of topological stacks \cite{No2} to manipulate Deligne-Mumford topological stacks as one would manipulate ordinary spaces.

These computational tools are used in the second section to revisit a well-known  class of examples of complex projective surfaces, namely the Hirzebruch 
covering surfaces branched along a line arrangement \cite{BHH}. The main original result in this section is:
\begin{theoi}
 Hirzebruch covering surfaces with equal weights do satisfy the Shafarevich conjecture. 
\end{theoi}
Hirzebruch covering surfaces  come indeed from the construction  given in the first section and the corresponding quasi-projective manifold is nothing but the 
complement of the union of the lines for a projective linear  arrangement - a rich fundamental group, indeed, and much studied. Of course, it would have been much better
to find a counterexample to the Shafarevich conjecture among such classical examples of surfaces, but  one has to try harder. 

The concluding section is devoted to the formulation of several open questions which are not 
investigated to the best of my knowledge in the litterature. I believe that the higher dimensional analogue of Theorem 2 
is, by far, the most accessible one. The group theoretical consequence of the Shafarevich conjecture for orbifolds (Question \ref{mainconj}) is
the promised reformulation of the main open question in \cite{BK}. Finally, Theorem \ref{shaflinfinal} gives a reduction of the Shafarevich conjecture for $X$
to a group theoretical statement concerning the images in the fundamental group of $X$ of the fundamental groups of the connected compact complex analytic 
spaces over $X$ for which a proof seems to be out of reach although a counterexample might be more accessible. 

\section{Quasi-K\"ahler groups, compact K\"ahler orbifolds, uniformizability and developpability.}

\subsection{The fundamental group of a weighted DCN}\label{wdcn}

Let $M$ be a (Hausdorff second countable) complex analytic space and $D$ be an effective Cartier divisor. Let $r\in \N^*$. Then, one can construct $P \to M$ the principal 
$\C^*$-bundle attached to $\mathcal{O}_M(-D)$ and the tautological section $s_D \in H^0(M,\mathcal{O}_M(D))$ can be lifted 
to a holomorphic function $f_D: P \to \C$ satisfying $f_D(p. \lambda)= \lambda f_D(p)$ for every $\lambda\in\C^*, \ p\in P$. 
Define a complex analytic space $Y:=Y_D\subset P \times \C_z$ by the equation $z^r=f_D(p)$. One can define a $\C^*$-action on $Y$ by
$(p,z). \lambda= (p\lambda^r, \lambda.z)$. The complex analytic stack $$M(\sqrt[r]{D}):=[Y_D/\C^*]$$ (see e.g. \cite{EySa, bn}) 
is then a Deligne-Mumford separated complex analytic stack with trivial generic isotropy groups whose moduli space is $M$
and an orbifold if $M$
and $D$ are smooth. The nontrivial isotropy groups live over the points of $D$ and are isomorphic to $\mu_r$
the group of $r$-roots of unity. In the smooth case, the corresponding differentiable stack can be expressed as the quotient by the natural infinitesimally free $U(1)$ action 
on the restriction $UY$ of $Y$ to the unit subbundle for (any hermitian metric) of $P$  which is a manifold indeed. It is straightforward to
see that this is an analytic version of Vistoli's root stack construction:
\begin{lem}
 If $(M, D)$ is the analytification of $(\mathcal{M}, \mathcal{D})$ a pair consisting of a $\C$-(separated) scheme and a Cartier divisor, 
 then $M(\sqrt[r]{D})$ is the analytification of $\mathcal{M}_{O (\mathcal{D}), s_{\mathcal{D}}, r}$ in the notation of \cite[Section 2]{cadman2007}. 
\end{lem}
 The main property I will use is treated in the scheme-theoretic setting by \cite{cadman2007}:  
 \begin{lem} If $S$ is a complex analytic stack 
 $$Hom(S,M(\sqrt[r]{D})) = \{f: S\to M \ \mathrm{holomorphic}, D_S \  \mathrm{Cartier \ on \ } S \  \ \mathrm{s. t.} \ D_S=r.f^*D\}. $$
 \end{lem}

Let $\bar X$ be a compact K\"ahler manifold and $x_0\in X$ a base point, $n=\dim_{\C}(X)$, 
and let $D:=D_1+ \ldots + D_l$ be a simple normal crossing divisor whose smooth irreducible components are denoted by $D_i$. We assume for simplicity $x_0\not \in D$. 
For each choice of weights $d:=(d_1, \ldots, d_l)$, $d_i\in \N^*$,  one may construct 
as \cite[Definition 2.2.4]{cadman2007} does in the setting of scheme theory 
the compact K\"ahler orbifold (Compact K\"ahler DM stack with trivial generic isotropy) 
$$\X(\bar X, D, d):= \bar X(\sqrt[d_1]{D_1} )\times_{\bar X} \ldots \times_{\bar X} \bar X(\sqrt[d_l]{D_l}).$$
In other words, $\X(\bar X, D, d)=[Y_{D_1}\times_{\bar X} \dots \times_{\bar X} Y_{D_l} / \C^{*l}]$. 
Denote by $X$ the quasi-K\"ahler manifold $X:=\bar X \setminus D$. View $\X(\bar X, D, d)$ as an orbifold compactification of $X$ and denote by $j_d: X \hookrightarrow \X(\bar X, D, d)$ the natural open immersion. 

By Zariski-Van Kampen, 
the fundamental group $\pi_1(\bar X,x_0)$ is the quotient of $\pi_1(X,x_0)$ by the normal subgroup generated by the $\gamma_i$, where $\gamma_i$ is a meridian loop for $D_i$. 

This generalizes to orbifolds, see e.g. \cite{No1}, and $\pi_1(\X(\bar X, D, d),x_0)$ is the quotient of $\pi_1(X,x_0)$ by the normal 
subgroup generated by the $\gamma_i^{d_i}$. Note that, if
$d=(1, \ldots, 1)$,  $\X(\bar X, D, d)=\bar X$. In particular the natural map $$j_{d*}:\pi_1(X,x_0)\to\pi_1( \X(\bar X, D, d),x_0)$$ is surjective. 

Say $d\in \N^{*l}$ divides $d'\in \N^{*l}$, which will be denoted by $d|d'$,  whenever $d'_i/d_i \in \N$ for all $i \in \N$. If $d|d'$, we have a map (1-morphism of DM stacks)
$p_{d,d'}: \X(\bar X, D, d')\to \X(\bar X, D, d)$ inducing the identity on the open substack $X$. This gives surjective group homomorphisms:
$$\pi_1(X) \buildrel{j_{d *}}\over \to \pi_1(\X(\bar X, D, d')) \buildrel{p_{d,d' *}}\over \to \pi_1(\X(\bar X, D, d))\buildrel{p_{1, d *}}\over \to \pi_1(\bar X).
$$

The same construction can be done for the \'etale fundamental group $\pi_1^{et}$ which is the profinite completion of $\pi_1$. 

The log pair $(\bar X, D)$ has to be pretty rigid if one wishes that the fundamental group $\pi_1(\X(X,D,d),x_0)$ be very different from $\pi_1(\bar X,x_0)$. Here is a simple manifestation: 
\begin{prop}\label{nori}
Assume $ \mathcal{O}_{D_q}(D_q)$ is ample for all $1\le q\le l$, and  $d|d'$. Then we have a finite central extension:
$$1 \to A \to \pi_1(\X(\bar X, D, d'),x_0)  \to \pi_1(\X(\bar X, D, d),x_0) \to 1
$$
by a finite abelian group having  generators $g_q$ whose order divides $d'_q /  d_q$, $q=1, \ldots,l$.
\end{prop}
\begin{prv}
Indeed $\gamma_1, \ldots, \gamma_l$ are central in $\pi_1(X)$ by \cite{Nori}.
\end{prv}

\subsection{A possible obstruction to residual finiteness for quasi-K\"ahler groups}

Consider the natural group homomorphism $$j_{\infty *}: \pi_1(X,x_0) \to \varprojlim (\pi_1(\X(\bar X, D, d),x_0), p_{d,d' *}).$$ 

\begin{lem} \label{rfinj}

One has $$\ker(j_{\infty*}) \subset \ker(\pi_1(X,x_0) \to \pi_1^{et}(X,x_0))= \bigcap_{ H \mathrm{finite \  index \  subgroup \ of \ } \pi_1(X,x_0)} H. $$ 
Hence, if $\pi_1(X,x)$ is residually finite, then
 $j_{\infty*}$ is injective.

\end{lem}

\begin{prv}  Let $\gamma\in \pi_1(X,x_0)\setminus \ker(\pi_1(X,x_0) \to \pi_1^{et}(X,x_0))$. Then there exists a finite group $G$ and a 
morphism $\phi: \pi_1(X,x_0)\to G$ such that $\phi(\gamma)\not =1$. Let $d_i$ be the order of $\phi(\gamma_i)$ and $d=(d_1, \ldots, d_q)$. 
Then $\phi(\gamma_i^{n_i})=1$ and $\phi$ factorizes as $\phi= \bar \phi\circ j_{d*}$ where $\bar \phi: \pi_1(\X(\bar X, D, d),x_0)\to G$. 
In particular $\bar \phi (j_{d*}(\gamma))=\phi(\gamma)\not =1$ and a fortiori $j_{d*}(\gamma)\not=1$ and also $j_{\infty*}(\gamma)\not =1$.

\end{prv}

 The proof of lemma \ref{rfinj} also gives:
\begin{coro}
 The natural morphism $$j_{\infty*}^{et}: \pi_1^{et}(X,x_0) \to \varprojlim (\pi^{et}_1(\X(\bar X, D, d),x_0), p_{d,d' *})$$
is an isomorphism. 
\end{coro}

The best statement in the reverse direction  I could think of is that if $j_{\infty*}$ is injective and $\pi_1(\X(\bar X, D, d),x_0)$ is residually finite for all $d$ sufficiently divisible then 
 $\pi_1(X,x_0)$ is residually finite.
 
\begin{exe}
The inclusion in  Lemma \ref{rfinj} may be strict. 
 \end{exe}
\begin{prv}
 Indeed let $Y$ be a compact K\"ahler manifold and $L$ be a holomorphic line bundle on $Y$. 
 Let $X$ be the total space of $L$ minus the zero section. One has a central extension:
 $$1 \to \Z \slash q \Z \to \pi_1(X) \to \pi_1(Y) \to 1. 
 $$
 Assume  that $q=0$ (e.g.: $\pi_2
(Y)=0$). One may write $X=\bar X \setminus D$ where $\bar {X}= \mathbb{P} (\C\oplus L)$ and $D=D_1 \cup D_2$ is the union of the zero and infinite section of $L$. 
One has $\pi_1(\X(\bar X, D, d))= \pi_1(X) \slash g.c.d(d_1, d_2) \Z$ and $j_{\infty*}$ 
is injective whereas $\pi_1(X,x_0)$ needs not be residually finite \cite{Del} \cite{Rag} . 
\end{prv}

\begin{lem} $\ker(j_{\infty*})\subset \pi_1(X,x_0)$ is independent of the log-smooth compactification of $X$. 
          
\end{lem}

\begin{ques}
Does there exists $(\bar X, D)$ such that $j_{\infty*}$ is not injective?
\end{ques}
\begin{exe} The  example in Proposition \ref{nonunif} satisfies   $\ker(j_{\infty*})=\{ 1 \}$. 
 
\end{exe}

\subsection{A developpable non-uniformizable compact K\"ahler orbifold}

An orbifold $\X$ is said to be {\em developpable} if its universal covering stack \cite{No2} \footnote{
The covering theory of \cite{No2} has nothing 
to do with the complex structure. In fact, it works for reasonable topological DM stacks (with finite isotropy) e.g.: DM stacks in the Complex Analytic category.}
is an ordinary manifold which is equivalent to the
injectivity of all local inertia morphism $I_x=\pi_1^{loc}(\X,x) \to \pi_1(\X)$, $x$ being an orbifold point of $\X$. It is said to be {\em uniformizable}
whenever the profinite completion $I_x \to \pi_1^{et}(\X)$  of all local inertia morphism is injective\footnote{The litterature also uses 
{\em good} orbifold for developpable and {\em very good} orbifold for uniformizable.}. When $\pi_1(\X)$ is residually finite, 
uniformizability and developpability are equivalent for $\X$.

The fundamental group of a compact K\"ahler uniformizable orbifold is K\"ahler, i.e.: 
occurs as the fundamental group of a compact K\"ahler manifold (see e.g.: \cite{vjm}). 
In loc. cit. it was asked for an example of 
 a developpable non-uniformizable K\"ahler orbifold. 

\begin{prop} \label{nonunif}
 There exists a pair $(\bar X, D)$ where $\bar X$ is complex projective, $D$ is a smooth divisor and $d\in \N^*$ such that $\X(\bar X,D,d)$ is developpable but not uniformizable. 
\end{prop}

\begin{prv} This is nothing but Toledo's famous example \cite{Tol}. 
Assume $n\ge 4$ and $n\equiv 0 [2]$. Assume $\bar X$ is uniformized by a symmetric domain $\Omega$ of type $IV$ and 
$D$ is a totally geodesic smooth complex hypersurface.

For every $d\in \N$,  
$\X(\bar X,D,d)$ is developpable. 
Indeed let $\pi^{u}: \Omega \to \bar X$ be the universal covering map. Then 
$\pi^{u \ -1}(D)$ is a countable disjoint union of totally geodesic complex hypersurfaces $(H_i)_{i\in \N}$. Those are equivalent
by the automorphism group of $\Omega$ to a hyperplane section $\{z_1=0\}$ of $\Omega$ embedded in $\C^n_{z_1, \ldots, z_n}$ by Harish-Chandra embedding.
Hence we may assume $H_0=\{z_1=0\}\cap \Omega$. Since $\Omega$ is convex in its Harish Chandra representation the pair $(\Omega, H_0)$ is diffeomorphic 
to $(\C, 0) \times \C^{n-1}$. 
The complex-analytic stack $\Omega\times_{\bar X} \X(\bar X,D,d)$ is, by functoriality of the root stack, equivalent to $\Omega[\sqrt[d]{H_0}]\times_{\Omega} \Omega[\sqrt[d]{H_1}]\times_{\Omega} \ldots$
and is a covering  of topological stacks of $\X(\bar X,D,d)$. Now consider $c_d: \C^n \to \C^n, \ (z_1,\ldots z_n)\mapsto (z_1^d, z_2, \ldots, z_n)$. Then, $c_d^{-1}(\Omega) \to \Omega [\sqrt[d]{H_0}]$ is a universal
covering stack and $\Omega [\sqrt[d]{H_0}]$ is uniformizable. It follows that $c_d^{-1}(\Omega) \times _{\Omega}  \Omega[\sqrt[d]{H_1}]\times_{\Omega} \ldots$ is a covering stack of $\X(\bar X,D,d)$ such that a connected component of the preimage 
of the substack supported on $D$ is a smooth hypersurface having no inertia groups. It follows that $\X(\bar X,D,d)$ is developpable.

We now have to prove that for every $d\in \N$ such that $d\not \in \{1,2,4,8\}$ 
$\X(\bar X,D,d)$ is not uniformizable. Choose a tubular neighborhood $V$ of $D$ in $\bar X$ which is diffeomorphic to the closed unit ball in the normal bundle 
of $D$ in $\bar X$ so that $D$ corresponds to the zero section. Let $N=\partial V$, $U=V\setminus D$. Since $\pi_2(D)=\{ 0 \}$, we have a central extension:
$$1\to \Z=\gamma^{\Z} \to \pi_1(N) \to \pi_1(D)\to 1$$
and $\pi_1(N)$ injects in $\pi_1(X)$ (see \cite{Tol}). 

Since $\pi_1(D)$ injects in $\pi_1(\bar X)$ and the inertia group at the base point $y$ of $D$ injects into $\pi_1(\X(\bar X,D,d))$, 
it follows that $\pi_1(N)/d\Z$ injects in $\pi_1(\X(\bar X,D,d))$ as the image of the fundamental group of the open substack $U[\sqrt[d]{D}]$. 
The central subgroup $\Z \slash d\Z\subset\pi_1(N)/d\Z$ identifies with the inertia group at the base point $y$ of $D$. 

The main fact used by Toledo is the following: every finite index subgroup in $\pi_1(N)$ contains $8\Z$  \cite{Rag}. Hence the kernel of any finite representation of $\pi_1(\X(\bar X,D,d))$
contains the subgroup $8\Z+d \Z \slash d \Z\subset \Z \slash d\Z= \pi_1^{loc}(\X(\bar X,D,d),y) \subset \pi_1(\X(\bar X,D,d))$. 
\end{prv}

\begin{ques}
 Are the orbifold K\"ahler groups $\pi_1(\X(\bar X,D,d))$ constructed above genuine K\"ahler groups for $d\not \in \{1,2,4,8\}$?
\end{ques}

\begin{rem}
Analysing the above argument, we see that $\X(\bar X,D,d)$ is developpable provided $(\pi^u)^* O_{\bar X}(D_i)$ is holomorphically trivial for all $1\le i \le l$ where  $\pi^u: \wt{\bar X^u} \to \bar X$ is
 the universal covering space of $\bar X$. This condition is satisfied e.g. when $ \wt{\bar X^u} $ is Stein and $[D_i]$ lies in the image of $H^2(\pi_1(\bar X), \Z) \to H^2(\bar X, \Z)$. 
 A sufficient condition for the latter is $\pi_2(\bar X)=\{0\}$.
\end{rem}

\subsection{Basic obstructions towards developpability}

It seems hopeless to classify uniformizable or developpable compact K\"ahler orbifolds. The $n=1$ case was known as the Fenchel-Nielsen problem and was solved by
Fox. 
 
The following trivial lemma enables to find obstructions for developpability or uniformizability.

\begin{lem}
 If $\X$ is a developpable (resp. uniformizable) compact K\"ahler orbifold, any connected closed substack is developpable (resp. uniformizable).
\end{lem}

Let $\gamma: \X \to \X^{mod}$ be  the canonical map to the moduli space of $\X$. Every closed substack of $\X$ has the form ${\mathcal Y}= Y\times_{\X^{mod}} \X$ where $Y$ is a closed substack of $\X^{mod}$. 
The only non-uniformizable (or non-developpable) compact complex orbi-curves being the orbifold $\X({\mathbb{P}}^1, \sqrt[d]{p})$, $d>1$ and 
$\X({\mathbb{P}}^1, \sqrt[d_1]{p_1})\times_{{\mathbb{P}}^1}\X({\mathbb{P}}^1, \sqrt[d_2]{p_2})$, $d_1\not = d_2$ both $ >1$, it follows that: 
\begin{coro}
 If $\bar X$ contains a smooth rational curve $C$ transverse to $D$ such that $C\cap D=\{p\}$ or $C\cap D_1= \{p_1\}$ $C\cap D_2=\{p_2\}$ with $d_1\not = d_2$,  both $ >1$, and $C\cap D_i=\emptyset$ if $k>2$, then $\X(\bar X,D,d)$ is not developpable. 
\end{coro}
One may formulate obstructions attached to singular or non transverse rational curves (say when the singular points are disjoint from $D$)  but I will refrain addressing this in complete generality.

The substack of $\X(X,D,d)$ corresponding to $D_p$ is equivalent to $B (\Z/d_p\Z)\times \mathcal{D}_p$ where $\mathcal{D}_p$ is the induced orbifold:
$$\mathcal{D}_p:=\X (D_p, \sum_{i=1, i\not=p}^l D_i\cap D_p, (d_1, \ldots, \hat d_p, \ldots)).$$

\begin{coro}
 If $\X(X,D,d)$ is uniformizable (resp. developpable), so is $\mathcal{D}_p$ for all $1\le p\le l$.
\end{coro}

There is no obstruction to developpability coming only from the fundamental group.

\begin{prop} \label{develo}
 Let $\X$ be a compact K\"ahler orbifold. Then there exists $\X'$ a developpable compact K\"ahler orbifold such
 that $\pi_1(\X')=\pi_1(\X)$. Furthermore the moduli space of $\X'$ has a proper holomorphic bimeromorphic mapping 
 to the moduli space of $\X$.
\end{prop}

\begin{prv}
 This is a corollary of Noohi's covering theory. Indeed, \cite{No2} constructs 
 a simply connected complex orbifold $\wt{\X^u}$ upon which $\Gamma=\pi_1(\X)$ acts with an equivalence $[\wt{\X^u}/\Gamma]=\X$. 
 We have the moduli map $\gamma:\wt{\X^u}\to \Xi:=(\wt{\X^u})^{mod}$. The complex space with quotient singularities $\Xi$ is simply connected since
 a covering of $\Xi$ lifts to a covering of $\wt{\X^u}$. Also, 
$\Gamma$ acts in properly dioscontinuous fashion on $\Xi$. The normal complex orbispace $\X^{red}:=[\Xi/\Gamma]$ has the same moduli space than $\X$ is sandwiched between $\X$ and $\X^{mod}$ and is the smallest such sandwiched 
 orbifold with fundamental group $\Gamma$. Its isotropy structure is a reduction of 
 the isotropy structure of $\X$, whence the terminology. 
 
 Consider $\Xi'\to \Xi$ a canonical (e.g.: Encinas-Villamayor) 
 resolution of the singularities of $\Xi$. Then, $\Xi'$ is simply connected
 by \cite[Sect. 7]{Kol1} \cite{Taka}. 
 Then $\X'=[\Xi'/\Gamma]$ is the required orbifold. 
 
 It should be noted that the K\"ahler form of $\X$ gives rise to a $\Gamma$-equivariant K\"ahler form on $\Xi$ 
 and then on $\Xi'$ using an orbifold variant of some standard techniques in M. P\u aun's thesis (see \cite{EySa}). Hence $\X'$ is K\"ahler indeed.
\end{prv}

\section{Hirzebruch covering surfaces}

Let us apply the previous construction to the study of the fundamental group of a famous class of algebraic surfaces introduced by Hirzebruch in \cite{Hirz}. 

\subsection{Reformulation of Hirzebruch's construction}
Let $\mathcal{A}=\{L_1, \ldots, L_k \}$ be a line arrangement in $\mathbb{P} ^2$. Let $\{ p_1, \ldots, p_s\}$ be the points of multiplicity $\ge 3$ in $\sum_{i} H_i$.
Denote by $\bar X$
the blow up at $p_1, \ldots, p_s$ of $\mathbb{P} ^2$. Denote by $D_1, \ldots, D_k$ the strict transforms of $L_1, \ldots, L_k$ and by $D_{k+1}, \ldots D_{k+s}$ 
the exceptional curves over $p_1, \ldots , p_s$.
Put $l=k+s$ and consider $D=D_1+ \ldots + D_l$. For $N\in \N$,  define $$\mathcal{M}_N(\mathcal{A}):= \X(\bar X, D, (N, \ldots, N)).$$ 

\begin{prop} \label{unifhirz} Unless $\mathcal{A}$ is a pencil,  $\mathcal{M}_N(\mathcal{A})$ is uniformizable. More precisely,
$H_1(X, \Z)\simeq \Z^{k-1}$ and the morphism $\pi_1(X) \to H_1(X, \Z/N\Z)=(\Z/N\Z)^{k-1}$
factors through a morphism $\eta:=\eta_{\mathcal{A}}: \pi_1(\mathcal{M}_N(\mathcal{A}))\to (\Z/N\Z)^{k-1}$ which is injective on the isotropy groups of $\mathcal{M}_N(\mathcal{A})$. 
\end{prop}
\begin{prv} The isomorphism $H_1(X, \Z)\simeq \Z^{k-1}$ is classical, see e.g. \cite{OT}. 
 
 Next,  consider the case $\mathcal{A}$ is an affine pencil, i.e.: a pencil $\mathcal{A}'$ plus a line which does not pass through the base point of the pencil.

 The case $k=3$ is somewhat simpler and I treat it first. In this case, $\mathcal{A}$ is a triangle, which can be described as $x_0x_1x_2=0$ where $(x_0,x_1,x_2)$ is a homogenous coordinate system of $\mathbb{P}^2$. The map 
 $\mathbb{P}^2 \to \mathbb{P}^2, \ [x_0:x_1:x_2]\to [x_0^N: x_1^N:x_2^N]$ lifts to an etale $(\Z/N\Z)^2$-covering of  $\mathcal{M}_N(\mathcal{A})$ 
 which is actually its universal covering space. The resulting isomorphism
 $\pi_1(\mathcal{M}_N(\mathcal{A}))\simeq (\Z/N\Z)^2$ is actually given by $\eta_{\mathcal{A}}$. 
 
 For the case $k\ge 4$,  one needs to blow up the base point of the pencil and $\bar{X}\simeq \mathbb{F}_1$ the  Hirzebruch surface of index $1$. Denote by
 $p:\mathbb{F}_1 \simeq \mathbb{P}(\mathcal{O}_{\mathbb{P}^1}\oplus\mathcal{O}_{\mathbb{P}^1}(1) )\to \mathbb{P}^1$. 
 Then $D$ is the union of $p^{-1}(S)$ with $S\subset \mathbb{P}^1$ is a finite subset of cardinal $k-1$, 
 the unique $-1$ curve $E$ and a $+1$ curve $F$.
 When $D'$ is a SNC divisor on $\mathbb{F}_1$ we will also denote by $p$ the composition $\X(\mathbb{F}_1, D',d) \to \mathbb{F}_1 \to \mathbb{P}^1$. 
 
 This implies that $\mathcal{M}_N(\mathcal{A}')$ is not uniformizable since the general fiber of $p$ is a teardrop. On the other hand, the natural map
 $p: \mathcal{X}(\bar X, p^{-1}(S), (N, \ldots, N))\to \mathbb{P}^1$ lifts to 
 $\bar p: \mathcal{X}(\bar X, p^{-1}(S), (N, \ldots, N))\to \mathbb{P}^1(\sqrt[N]{S})$ by \cite[Remark2.2.2]{cadman2007}. The map $\bar p$ 
is a $\mathbb{P}^1$-bundle and yields an isomorphism on $\pi_1$. Furthermore, $$H_1(\mathbb{P}^1(\sqrt[N]{S}) , \Z/ N. \Z)= (\Z/N\Z)^{k-2}$$
and the resulting map $\pi_1(\mathbb{P}^1(\sqrt[N]{S}) )\to (\Z/N\Z)^{k-2}$ is injective on isotropy groups so that there is 
a smooth abelian uniformization of degree $N^{k-2}$,  $e: C\to \mathbb{P}^1(\sqrt[N]{S}) )$. 
Then $C\times_{\mathbb{P}^1} \mathcal{M}_N(\mathcal{A})$ is $N$-th root stack of $Y=\mathbb{P}(\mathcal{O}_C \oplus  e^*\mathcal{O}_{\mathbb{P}^1}(1))$ along the divisor 
corresponding to the two tautological sections. Since $e^* \mathcal{O}_{\mathbb{P}^1}(1)$ has degree 
$N^{k-2}$ we may construct $M$ an invertible sheaf on $C$ such that $M^{\otimes N}\simeq e^*\mathcal{O}_{\mathbb{P}^1}(1)$. 
This is exactly what is required to construct a fiberwise $N$-th power map $\mathbb{P}(\mathcal{O}_C \oplus M) \to Y$ which lifts to an etale cover of $C\times_{\mathbb{P}^1}\mathcal{M}_N(\mathcal{A})$ compatible
with the projection to $C$. 

By composition, we get an etale covering map $e':\mathbb{P}(\mathcal{O}_C \oplus M) \to \mathcal{M}_N(\mathcal{A})$ in particular $\mathcal{M}_N(\mathcal{A})$ is uniformizable. 

It remains to be seen  that the etale covering map $e'$ corresponds precisely to $\ker(\eta_A)$. To see this observe that 
$p: \mathcal{M}_N(\mathcal{A}) \to \mathbb{P}^1(\sqrt[N]{S})$ is a fibration of DM topological stacks whose fiber is a football $\mathbb{P}^1(\sqrt[N]{0+\infty})$ and, thanks to \cite{No2} gives rise to an exact sequence: 
$$1\to \Z/ N \Z \to \mathcal{M}_N(\mathcal{A}) \to \pi_1(\mathbb{P}^1(\sqrt[N]{S}))\to 1.$$ 
This sequence splits since $p$ has a section and the monodromy action is trivial since $p$ is a 
pull back  of a fibration  over $\mathbb{P}^1$, which has a trivial monodromy action since $\mathbb{P}^1$ is simply connected. 
Hence, the exact sequence underlies a product decomposition $\pi_1(\mathcal{M}_N(\mathcal{A}))= \Z/ N\Z \times  \pi_1(\mathbb{P}^1(\sqrt[N]{S}))$ and the claim follows. 
 
 Now to the general case. For every double point $p$ of $\mathcal{A}$, since  $\mathcal{A}$ is not a pencil, there a triangle subarrangement $\mathcal{A'}$ of  $\mathcal{A}$
 and the map $\mathcal{M}_N(\mathcal{A}) \to \mathcal{M}_N(\mathcal{A}) $ is an isomorphism near $p$. Since $\eta_{\mathcal{A}'}$ is a quotient of  $\eta_{\mathcal{A}}$ it follows from the triangle case that
 $\eta_{\mathcal{A}}$ is injective on the isotropy group at $p$. For every $m\ge 3$-uple point $q$, there is a near pencil $\mathcal{A}'$ subarrangement with the $m$ lines of $\mathcal{A}$ through $p$ and the 
 statement follows from the case of $\mathcal{A}'$ by the same argument. 
 \end{prv}

 \begin{defi}
  The covering surface $M_N(\mathcal{A}) \to \mathcal{M}_N(\mathcal{A})$ corresponding to $\ker (\eta_{\mathcal{A}})$  is a smooth projective surface with an action of $G:=(\Z/N\Z)^{k-1}$
  such that $\mathcal{M}_N(\mathcal{A})\simeq [M_N(\mathcal{A}) \slash G]$ and is called the Hirzebruch covering surface branched with index $N$ over $\mathcal{A}$. 
  \end{defi}

The surface $M_N(\mathcal{A})$ is the \lq Kummersche \"Uberlagerung der projektiven Ebene\rq \ described  in \cite[Kapitel 3]{BHH}. 

\subsection{The fundamental group of Hirzebruch covering surfaces} \label{hirzshaf} This reference \cite{BHH} contains a wealth of information on $M_N(\mathcal{A})$.
However, the fundamental group of  $M_N(\mathcal{A})$ does not seem to have been studied systematically in the litterature.
One nevertheless knows quite a lot. 

The computation of $b_1(M_N(\mathcal{A}))$ in terms of the characteristic variety of $\mathcal{A}$ was carried out by \cite{Hi,Sa} see also \cite{Su}.
It  may be possible to study the LCS and Chen ranks of $\pi_1(M_N(\mathcal{A}))$ along the lines of \cite{Su}. 

Here, I will give a necessary and sufficient condition for the fundamental group of  $M_N(\mathcal{A})$ to be finite.
However, I will begin by describing the main examples.

When $\mathcal{A}$ is a general type arrangement $\pi_1(M_N(\mathcal{A}))$ is finite abelian by Corollary \ref{nori}.

If $p$ is a $k\ge 3$-uple point of $\mathcal{A}$
the previous construction gives a surjective map $f_p:\mathcal{M}_N(\mathcal{A}) \to \mathbb{P}^1(\sqrt[N]{S_p})$ where $\# S_p=k$ with connected fibers hence the induced mapping on $\pi_1$ is surjective (it actually has a section). 
It follows that if $N\ge 3$ or $N=2$ and $k\ge 4$ then $\pi_1(M_N(\mathcal{A}))$ is infinite and has an infinite image linear representation in the presence of singular points. 

\begin{exe}
 Let $\mathcal{B}^* \subset \C^4$ be the braid arrangement whose equation is $$\mathcal{B}^*=\{(x_1,\ldots, x_4) | \  \prod_{1\le i<j\le 4}(x_i-x_j)=0\}.$$ The intersection of these hyperplanes is the line $\C. (1,1,1,1)$.  Let $\mathcal{B}$ be the line arrangement 
 induced by the central arrangement  $\mathcal{B}^*/ \C.(1,1,1,1)$ in $\C^4/ \C.(1,1,1,1)\simeq \C^3$. Then $\mathcal{B}$ is the complete quadrilateral formed by the 6 lines joining every 2 among 4 points in general position. 
 One of the first computations in \cite{BHH} shows that $M_2(\mathcal{B})$ is a K3 surface (their computation gives $K=0, c_2=24$), hence $\# \pi_1(\mathcal{M}_2(\mathcal{B}))=2^5=32$. 
 \end{exe}

\begin{theo}\label{mainq}
 If $\mathcal{A}$ is not a pencil and has only double and triple points,  $\pi_1(\mathcal{M}_2(\mathcal{A}))\simeq (\Z/2\Z)^{\# \mathcal{A} -1} $. Consequently, $M_2(\mathcal{A})$ is simply connected. 
\end{theo}

The following appears in \cite{Lef1} with no proof except a reference to a personal communication I made to the author:  

\begin{coro} 
 The fundamental group $\pi_1(\mathcal{M}_2(\mathcal{A}))$ is finite 
 if and only if $\mathcal{A}$ has only double points or $\mathcal{A}$ has only double and triple points and $N=2$. 
\end{coro}

\begin{prv} 
Let $B_1, \ldots, B_r$ be small disjoint closed coordinate balls centered on the triple points of $\mathcal{A}$ and let $S_1, \ldots, S_r$ denote their boundary spheres. 
Define also $\Omega={\mathbb P}^2 \setminus (B_1\cup \ldots \cup B_r)$,  $\gamma: \mathcal{M}_2(\mathcal{A}) \to \mathbb{P}^2$ the moduli map, and 
$\Omega_2=\gamma^{-1}(\Omega)$, $\Sigma_{2,i}=\gamma^{-1} (S_i)$, $B_{2,i}=\gamma^{-1}(B_i)$ all viewed as topological DM stacks. 
By the Seifert-Van Kampen for orbifolds \cite{Hae}  \cite[p. 140]{Kap} (one may also follow the method of \cite{Zoo}) :
$$
\pi_1(\mathcal{M}_2(\mathcal{A}))=\pi_1(\Omega_2) *_{\pi_1(S_{2,1})} \pi_1(B_{2,1}) *_{\pi_1(S_{2,2})} \ldots *_{\pi_1(S_{2,r})} \pi_1(B_{2,r}). 
$$

On the other hand, let $\mathcal{A}'$ be a small perturbation of $\mathcal{A}$ so that $\mathcal{A}'$ is a generic arrangement. 
Define as above $\Omega={\mathbb P}^2 \setminus (B_1\cup \ldots \cup B_r)$,  $\gamma': \mathcal{M}_2(\mathcal{A}') \to \mathbb{P}^2$ the moduli map, and 
$\Omega'_2=(\gamma')^{-1}(\Omega)$, $\Sigma'_{2,i}=(\gamma')^{-1} (S_i)$, $B_{2,i}'=(\gamma')^{-1}(B_i')$ all viewed as topological DM stacks. 
By Seifert-Van Kampen:
$$
\pi_1(\mathcal{M}_2(\mathcal{A}'))=\pi_1(\Omega'_2) *_{\pi_1(S'_{2,1})} \pi_1(B'_{2,1}) *_{\pi_1(S'_{2,2})} \ldots *_{\pi_1(S'_{2,r})} \pi_1(B'_{2,r}). 
$$
By  Ehresmann lemma, we obtain a homeomorphism:
$$
(\Omega_2, \Sigma_{2,1}, \ldots, \Sigma_{2,r}) \simeq (\Omega'_2, \Sigma'_{2,1}, \ldots, \Sigma'_{2,r}).
$$
Hence, it is enough to show that the group morphisms $\pi_1(S_{2,i})\to \pi_1(B_{2,i})$ and $\pi_1(S'_{2,i})\to \pi_1(B'_{2,i})$
are isomorphic. Since $\pi_1(B'_{2,i})=H_1(B'_{2,i})=(\Z \slash 2 \Z)^3$ it is enough to show:
\begin{lem}
  $\pi_1(B_{2,i})=H_1(B_{2,i})=(\Z \slash 2 \Z)^3$. 
\end{lem}
\begin{prv} The proof of Proposition \ref{unifhirz} also yields that $H_1(B_{2,i})=(\Z \slash 2 \Z)^3$. 
Let $S\subset {\mathbb P}^1$ be a finite subset with $3$ elements. Then $B_{2,i}$ is a 
smooth fiber bundle over ${\mathbb P}^{1}(\sqrt[2]{S})$ with fiber $\Delta(\sqrt[2]{0})$ where $\Delta$ is an open unit disk.  
In particular by \cite{No2} we get an exact sequence:
$$ \pi_1(\Delta(\sqrt[2]{0}))=\Z \slash 2 \Z \to \pi_1(B_{2,i}) \to \pi_1({\mathbb P}^{1}(\sqrt[2]{S})) \to \{1\}.
$$
Now $\pi_1({\mathbb P}^{1}(\sqrt[2]{S}))=<a,b,c| a^2=b^2=c^2=abc=1>$ is easily seen to be abelian isomorphic to $(\Z\slash 2\Z)^2$. Hence, 
$\pi_1(B_{2,i}) $ is finite of cardinal 4 or 8. Since its abelianization has cardinal 8, it follows that 
$\pi_1(B_{2,i}) $ has cardinal 8, hence is isomorphic to its abelianization. 
\end{prv}
The proof of Theorem
\ref{mainq} is complete. 
\end{prv}
 
 \subsection{The Shafarevich conjecture for $M_N(\mathcal{A})$}
 
 \begin{theo} \label{shafhirzsurf} If $\mathcal{A}$ is not a pencil, the universal covering space of $M_N(\mathcal{A})$ is holomorphically convex.
 \end{theo}
 
 \begin{prv}

 We will first treat the $N\ge 3$ case. 
 The idea of the proof stems from Sommese's construction of minimal general type  surfaces with a Miyaoka-Yau ratio $c_1^2/c_2$ close to $3$ \cite{Som}.
 
 We may assume $\mathcal {A}$ is not generic  since the fundamental group is then finite and the universal covering space compact hence proper over the point which is Stein. 
 Let $\{p_1, \ldots, p_s\}$ be the  points of multiplicity $\ge 3$ in $\mathcal{A}$ and consider
 $F: \prod_{m=1}^s f_{p_m}: \mathcal{M}_N( \mathcal{A})\to \prod_{m=1}^s \mathbb{P}^1(\sqrt[N]{S_{p_m}})$. 
 
 The first case to be treated is if there are 3 non collinear such  points in $\mathcal{A}$.  By construction, $F$ contracts no curve. In particular $F$ is finite over an orbifold 
 whose universal covering space is a product of disks and complex lines and the corresponding covering space of $\mathcal{M}_N( \mathcal{A})$ is Stein. Hence, its universal covering space which is the same as the
 universal covering space of $M_N(\mathcal{A})$ is Stein too.
 
 Then we have to look at the case when all multiple points lie on a line $L$ that may or may not belong to $\mathcal{A}$.
 
A degenerate case is when there is a single  point of multiplicity $\ge 3$. In this case, since all 
the lines that do not pass through this multiple point are transverse to the arrangement, we may
delete them  thanks to lemma \ref{nori}. Indeed, if $f:\mathcal{Y}\to \mathcal{Z}$ is a {\em finite}  map of compact K\"ahler orbifolds such that the mapping induced by $f$ on $\pi_1$ has finite kernel
then the map from the universal cover of $\mathcal{Y}$ to the universal cover of $\mathcal{Z}$ is finite and  holomorphic convexity of the universal cover of $\mathcal{Z}$  implies
 holomorphic convexity of the universal cover of $\mathcal{Y}$.
This reduces us to the case of a near pencil where the analysis in  the proof of proposition \ref{unifhirz} shows that $M_N(\mathcal{A})$ is a ruled surface, the ruling being essentially given by $F$. 
In this case, the $\Gamma$-dimension of $M_N(\mathcal{A})$ is $1$. 

Now assume there are at least $2$ points of multiplicity $\ge 3$ all lying on $L$.
Once again, we may delete from the arrangement all the lines that do not pass through these points since they are transverse to the rest of the arrangement. 
We are reduced to the case of an arrangement which is the union of $s$ pencils $\mathcal{P}_1, \ldots, \mathcal{P}_s$ with at least $3$ lines. 
Actually, there are two subcases whether  the line $L$ does belong to these pencils or does not. 
The base points of the $\mathcal{P}_i$ coincide with the $p_i$ and do  belong to $L$. 

The proper transform $L'$ of $L$ in $\bar X=Bl_{p_1, \ldots,p_s}({\mathbb P}^2)$ is the only curve that is contracted by the map $F^{mod}$ induced by  $F$ on the moduli spaces.
In fact, $L'$ is the only positive dimensional fiber of $F^{mod}$. Denote by $\mathcal{L}'$ the closed substack of $\mathcal{M}_N( \mathcal{A})$ 
defined as the preimage of $L'$ by the moduli map.

\begin{lem} \label{finiteimage}
  $\mathrm{Im}(\pi_1(\mathcal{L}')\to \pi_1(\mathcal{M}_N( \mathcal{A})))$
is finite. 
\end{lem}

\begin{prv}
In the present case it is possible to compute the fundamental group of $\mathcal{M}_N(\mathcal{A})$. 

First  consider $X'=\bar X -\psi^{-1} L$ where $\psi:\bar X \to \mathbb{P}^2$ is the natural map. $X'\simeq \C^2$ has a SNC divisor $D'$ consisting of the union 
of $m$ parallel arrangements $P'_1, \ldots, P'_m$ in general position with respect to each other. 
In particular, by \cite{Fan},  $$\pi_1(X'\setminus D')\simeq \prod_i \pi_1(X'\setminus P'_i),$$ each factor being a free group on 
$\# P'_i$ generators. 

We have a $\pi_1$-surjective mapping $j:\X(X',D',(N, \ldots, N))\to \mathcal{M}_N(\mathcal{A})$.
The induced morphism of fundamental groups gives rise to  a surjective group morphism $\psi: \prod_i \pi_1 (\X(X', P'_i, N))\to \pi_1(\mathcal{M}_N(\mathcal{A}))$ and we denote by $\psi_i: \pi_1 (\X(X', P'_i, N))\to \pi_1(\mathcal{M}_N(\mathcal{A}))$ the $i$-th component. 
Here the group $ \pi_1 (\X(X', P'_i, N))$ is a free product of $\#P'_i$ copies of $\Z/N\Z$. Call the corresponding generators $\gamma_{i1}, \ldots, \gamma_{i\# P_i}$ since they correspond to the meridian loops of the lines in $P_i$ . 

There is also  a group morphism $\pi_1(f_{p_i}): \pi_1(\mathcal{M}_N(\mathcal{A}))\to  \pi_1(\mathbb{P}^1(\sqrt[N]{S_{p_i}}))$. 
The group morphism $\pi_1(f_{p_i}) \circ \psi_i: \pi_1 (\X(X', P'_i, N)) \to \pi_1(\mathbb{P}^1(\sqrt[N]{S_{p_i}}))$ can be decribed as the quotient map 
by the normal subgroup  $H_i$ generated by $\gamma_{i1}.\gamma_{i2}. \ldots \gamma_{i\# P_i}$ if $L$ is not in the pencil or the normal
subgroup generated by $(\gamma_{i1}.\gamma_{i2}. \ldots \gamma_{i\# P_i})^{N}$ if so. 

In particular  $\pi_1(F)\circ \psi$ is surjective and $\ker(\pi_1(F)\circ \psi)=\prod_i H_i$. It follows that $\pi_1(F)$ is surjective and that the kernel $H$ of $\psi$ is contained 
in
$\prod_i H_i$ so that $\psi$ presents $\pi_1(\mathcal{M}_N(\mathcal{A}))$ as the quotient group by $H$. But it is clear from the description above that $\prod_i H_i \subset \ker(\psi)$.
Hence $\pi_1(F): \pi_1(\mathcal{M}_N(\mathcal{A}))\to \prod_i \pi_1(\mathbb{P}^1(\sqrt[N]{S_{p_i}}))$ is an isomorphism. 
'

The substack  $\mathcal{L}'\subset\mathcal{M}_N(\mathcal{A})$ corresponding to $L'$ maps via $F$ to  an orbifold point of $\prod_i \mathbb{P}^1(\sqrt[N]{S_{p_i}})$ in particular
the image of $\pi_1(\mathcal{L}')$ by $\pi_1(F)$ is finite. Since $\pi_1(F)$ is an isomorphism, 
$\mathrm{Im}(\pi_1(\mathcal{L}')\to \pi_1(\mathcal{M}_N( \mathcal{A})))$ 
is finite, as claimed.
\end{prv}

It follows that the preimage of $L'$ in the universal covering space $\pi; \wt{M_N(\mathcal{A})^u}\to \mathcal{M}_N(\mathcal{A})$
is a union of compact connected components (that are easily seen to be smooth). 
The fibers of $F\circ \pi$ are actually either discrete or preimages of $L'$ in particular, these fibers have all 
their connected components compact. This enables to construct a proper holomorphic mapping
$s:\wt{M_N(\mathcal{A})^u} \to S:= S(\wt{M_N(\mathcal{A})^u} )$ contracting precisely the compact connected components of the preimage of $L'$.

The action of $\pi:=\pi_1({\mathcal M}_N(\mathcal{A}))$ descends to $S$ an commutes with $s$
and we have a finite map of complex analytic stacks $[S/\pi] \to \prod_{m=1}^s \mathbb{P}^1(\sqrt[N]{S_{p_m}})$. Since the latter stack is a product of hyperbolic and elliptic orbifolds it follows that there is an intermediate covering
$S\to S'$ with $S'$ finite over a product of disks and complex lines. Such an $S'$ is Stein and since a covering of a Stein space is Stein, $S$ itself is Stein. 

The case when $N=2$ can be treated as follows. 
First of all,  theorem \ref{mainq}   reduces the question to the case when there are $\ge 4$-uple points. 
The case where they are not all collinear or where there is at most one such point is  treated as  
 above. The case when they are collinear, say lie on some line $L$,  and
 there are some triple points outside $L$  can be treated using the idea of theorem  \ref{mainq} by moving slightly one line of a pencil
based at infinity to eliminate these triple points (this will not change the fundamental groups).  
Then one can apply \cite{Fan} as above to conclude.

\end{prv}

\begin{rem}
 Even if the group in Lemma \ref{finiteimage} was infinite, applying \cite{Nap} would give the conclusion. 
\end{rem}

 \section{Final remarks and open questions}
 \subsection{Group theoretical nature of $\pi_1({\mathcal M}_N(\mathcal{A}))$}
 Let us recall a well-known interesting  open question on arrangement groups:
 \begin{ques}  Given $\mathcal{A}$ a line arrangement, 
  is $\pi_1(X)=\pi_1(\mathbb{P}^2\setminus \cup_{H \in \mathcal{A}} H)$ residually finite? A linear group? 
 \end{ques}
The standard notation is $X(\mathcal{A}):=\mathbb{P}^2\setminus \cup_{H \in \mathcal{A}} H$. 

\begin{exe} There is an elementary isomorphism
 $\pi_1(X(\mathcal{B}))\simeq \bar P_4$ where $P_4=\Z \times \bar P_4$ is the decomposition of Artin's pure braid group $\pi_1(\C^4\setminus \mathcal{B}^*)$ as the product of its infinite cyclic center $\mathcal{Z}\simeq \Z$
 by the quotient group $\bar P_4:=P_4/\mathcal{Z}$. Since $B_4$ is linear, it follows that $\bar P_4$ is linear too. 
\end{exe}

A similar statement holds for the projectivisation $\mathcal{B}_n$ of the generic $3$-dimensional linear section of the (reduced) higher braid arrangements $\mathcal{B}_n^*$. 
The fundamental group of $\C^n \setminus{\mathcal{B}}^*_n$ is Artin's pure braid group $P_n$, $n\ge 4$. 

In a similar fashion, we ask:

\begin{ques} Given $N\in \N^*$ and $\mathcal{A}$ an arrangement, is $\pi_1({\mathcal M}_N( \mathcal{A} ))$ residually finite? A linear group? 
 \end{ques}
 
 In general $\pi_1(F): \pi_1(\mathcal{M}_N(\mathcal{A} ))\to \prod_i \pi_1(C_i)$ is not an isomorphism onto a product of orbicurves groups. Hirzebruch's famous example of 
a compact complex hyperbolic orbifold, $\mathcal{M}_5(\mathcal{B})$, gives a counterexample. Indeed the rational cohomological dimension of $\pi_1(\mathcal{M}_5(\mathcal{B}))$
is $4$ whereas, in this case, the rational cohomological dimension of $\prod_i \pi_1(C_i)$ is $8$! One may hope that the projection to 2 factors be an isomorphism but Zuo's 
computation of $b_1$ easily allows to exclude this possibility. 

Furthermore, the components of $F$ need no be the only orbicurve
fiberings of $\mathcal{M}_N(\mathcal{A})$.  There is indeed a fifth such fibering for $\mathcal{M}_N(\mathcal{B})$ . An easy way to see it is that in this case $\bar X$ is isomorphic to
$\overline{\mathcal{M}_{0,5}}$ the Kontsevich moduli space of $5$-pointed genus $0$ curves, hence $\mathcal{M}_N(\mathcal{B})$ has $\mathfrak{S}_5$-symmetry. There is $\mathfrak{S}_4$ subgroup
fixing the exceptional curve and lifting the $\mathfrak{S}_4$- symmetry of $\mathcal{B}$. Precomposing one of this obvious fiberings with a generator of the cyclic subgroup gives the fifth fibering. 

I did not investigate whether $\mathcal{M}_N(\mathcal{B})$ is a linear (resp. residually finite) group if $N\not=\{2, 5\}$. 

\subsection{Higher dimensional analogue of Theorem \ref{shafhirzsurf}}
For every projective hyperplane arrangement
$\cup_{H\in \mathcal{A}} H \subset \mathbb{P}^n$ we have the De Concini-Procesi compatification $\bar X \to \mathbb{P}^n$ obtained by blowing up certain
non-generic intersections of the hyperplanes in $\mathcal{A}$. The preimage $D$ of  $\cup_{H\in \mathcal{A}} H $ in $\bar X$ is a simple normal crossing divisor
whose components are labelled by the above non-generic intersections of hyperplanes. 
If $\mathcal{A}$ contains a subarrangement consisting in $n+1$ coordinate hyperplanes Proposition \ref{unifhirz} holds
to the effect that ${\mathcal M}_N({\mathcal A}):=\mathcal{X}(\bar X, D, (N, \ldots, N))$ is uniformizable.
The Zariski-Lefschetz theorem also reduces the finiteness of $\pi_1({\mathcal M}_N({\mathcal A}))$ to the case $n=2$. However,
I don't have a complete argument for  the Shafarevich conjecture in that case:

\begin{conj}
 Is the universal covering space of ${\mathcal M}_N({\mathcal A})$ holomorphically convex? 
\end{conj}
 I do believe an affirmative answer should not be difficult to get. 

 \subsection{A group theoretical consequence of the Shafarevich conjecture}

 Let $(\bar X, D,d)$ be as in subsection   \ref{wdcn}. For $I=\{ i_1, \ldots, i_r \} \subset \{ 1, \ldots, l \}$ we denote by $D_I$ the intersection $D_{i_1}\cap  \dots\cap D_{i_r}$. 
Since $D_I$ need not be irreducible we denote by $D_{Ii}$ $1\le i \le \alpha_i$ (where $\alpha_i=0$ is actually allowed!) its irreducible components. 

\begin{defi}
A positive dimensional $D_{Ii}$ is said to be exceptional if $$\mathrm{Im}(\pi_1(\mathcal{D}_{Ii}) \to \pi_1(\X(\bar X, D, d))$$ is a finite group, where its preimage in $\X(\bar X, D, d)$ is denoted by $\mathcal{D}_{Ii}$.
If $D_{Ii}$ is exceptional we shall say that $\mathcal{D}_{Ii}$ is an exceptional component. 
\end{defi}

The following is a consequence of the Shafarevich conjecture generalized to K\"ahler orbifolds. The case when $\bar X$ is projective and $\X(\bar X, D, d)$ is uniformizable
would be a consequence of the usual statement
of the conjecture for complex projective manifolds. 

\begin{ques}\label{mainconj}
 Let $\mathcal{D}'\subset \mathcal {D}$ be a connected union of exceptional components. Then $\Gamma':=\mathrm{Im}(\pi_1(\mathcal{D}') \to \pi_1(\X(\bar X, D, d)))$ is a finite group.
\end{ques}

It follows from \cite{EKPR, CCE} that for all $\rho: \pi_1(\X(\bar X, D, d)))\to GL_N(\C)$,  $\rho(\Gamma')$ is finite.

The case where $(\bar X, D)$ comes from an arrangement of lines with unequal weights 
does not seem to be a direct consequence of the
litterature on fundamental groups of complements of line arrangements. 
More generally, any sufficiently complicated configuration of curves in a smooth surface - meaning that the fundamental 
group of the complement has exponential growth with respect to some system of generators- will give rise to a non-trivial instance of 
this question. 

\subsection{Connection with the Bogomolov-Katzarkov examples}
The situation considered by \cite{BK} can be described as follows. Their $\bar X$ is a smooth algebraic surface with a fibration  $\pi:\bar X \to S$,  $S$ being a complete curve,
whose singular (scheme-theoretic) fibers are stable curves of genus $g\ge 2$. 
This data is equivalent to a mapping $S\to \overline{\mathcal{M}_g}$ which cuts nicely the boundary components
of $\overline{\mathcal{M}_g}$, the Deligne-Mumford moduli stack of stable curves of genus $g$. 
Their $D$ is the preimage of the critical set of $\pi$. And they set $d=(p, \ldots, p)$ where $p$
is a large prime number (e.g.: $p\ge 641$) so that the Burnside group on $2$ (or more) generators with exponent $p$ is infinite. 
They obtain that
 $\X=\X(\bar X, D, p)$ is uniformizable.
Let $B\subset S$ be the set of critical points of $\pi$ and denote by $\mathcal{S}$ the orbicurve $\mathcal{S}=\X(S, \sqrt[p]{B})$.
They prove that the surjective morphism $\pi_1(\X)\to \pi_1(\mathcal{S})$ may have an infinite kernel.
The open problem they pose is nothing but an instance of Question \ref{mainconj}.

\subsection{A strenghtened form of the Shafarevich Conjecture on Holomorphic Convexity}

Let $\Gamma$ be a finitely presented group and $G/ \Q$ be an affine algebraic group. Then there exists a finite type affine $\Q$-scheme ${\mathrm{Rep}}(\Gamma, G)$
which represents the functor 
$$
\left(\begin{array}{rcc}
 \Q-{\mathbf{algebras}} & \to & {\mathbf{Sets}} \\
 A & \mapsto & {\mathbf{Hom_{Grp}}} (\Gamma, G(A))
\end{array}\right).
$$
Here, an algebra is assumed to be commutative and associative. Define $A_{\Gamma, G}:=\Q[{\mathrm{Rep}}(\Gamma, G)]$. By the very definition 
of a representing functor there exists a universal representation $\rho_{\Gamma, G}: \Gamma \to G(A_{\Gamma, G})$ such that
every group morphism $\Gamma \to G(B)$ factors as $G(\phi)\circ \rho_{\Gamma, G}$ for some algebra morphism $\phi:A_{\Gamma, G} \to B$
where we identify $G$ and its  functor of points. 

\begin{defi} $\rho_{\Gamma, G}$ will be called
 the universal representation of $\Gamma$ with values in $G$. Also, for $N\in \N^*$,  $\rho_{\Gamma, N}:=\rho_{\Gamma, GL_N}$. 
\end{defi}

\begin{defi} Let $X$ be a connected quasi-K\"ahler manifold\footnote{Quasi-K\"ahler means realizable as a Zariski open subset of a compact K\"ahler manifold. }. Define Property SSC(X) holds as follows: 

For every $f:Z\to X$ where $Z$ is a compact connected complex analytic space
and $f$ is holomorphic $$\#\mathrm{Im}(\pi_1(Z) \to \pi_1(X)) <+\infty \Leftrightarrow
\forall N \in \N^* \ \#\rho_{\pi_1(X), N}(\mathrm{Im}(\pi_1(Z) \to \pi_1(X))<+\infty. $$
 \end{defi}

\begin{theo}\label{shaflinfinal} Let $X$ be a compact K\"ahler manifold such that SSC(X) holds. Then the universal covering manifold of $X$ is holomorphically convex. 
 \end{theo}
 \begin{prv} Denote by $\pi: \xu \to X$ the universal covering map. 
  By \cite{CCE} (see also \cite{EKPR} in the projective case) $\wtx_N=\xu / \ker(\rho_{\pi_1(X), N})$ is holomorphically convex. 
  Let $\wtx_N \to S_N$ be its Cartan-Remmert reduction. The group $\Gamma_N=\pi_1(X) / \ker(\rho_{\pi_1(X), N})$ acts in a proper discontinuous cocompact fashion on $S_N$. 
  Following \cite{Cam1,Kol1}, we consider the natural holomorphic mapping  $sh_N: X \to Sh_N(X):=S_N / \Gamma_N$, which is a Shafarevich morphism (aka.  $\Gamma$-reduction). 
  
  By the characteristic property of the Shafarevich morphism \cite{pasy}, there is a tower of morphisms of normal complex-analytic spaces under $X$
  $$ \ldots \to Sh_{N+1}(X) \to Sh_{N}(X) \to \ldots \to Sh_1(X)
  $$
  and by the ascending chain condition for complex-analytic subspaces of $X\times X$ applied to $\{X\times_{Sh_{N}(X)} X\}_{N\in \N^*}$ there is $N_0\in \N^*$
  such that for all $N\ge N_0$ $sh_N(X)=sh_{N_0}(X)$. 
  \begin{defi}
 $Sh_{\infty}(X):=sh_{N_0}(X)$ and $(sh_{\infty}: X \to Sh_{\infty}(X)):=sh_{N_0}$ will be called the linear Shafarevich morphism of $X$. 
  \end{defi}
  
  \begin{lem} Under
   SSC(X),  every fiber of $q:=sh_{\infty} \circ \pi$ has no non-compact connected component. 
  \end{lem}
  
  \begin{prv}
  Every fiber of $q$ is of the form $\pi^{-1}(Z)$ where $Z\subset X$ is a fiber of $sh_N$ for $N\ge N_0$ and is contracted by $sh_N$ for all $N$. By the 
  the characteristic property of $sh_N$,  $\rho_{\pi_1(X), N} (\mathrm{Im}(\pi_1(Z) \to \pi_1(X)))$ is finite for all $N$. Hence by SSC(X) $\mathrm{Im}(\pi_1(Z) \to \pi_1(X))$ is finite
  too, so that all connected components of the lift of $Z$ to $\xu$ are compact. 
 \end{prv}
 
 Hence, there is a proper morphism $s:\xu \to S_{\infty}$ where $S_{\infty}$ is a normal complex-analytic space whose points are in 1-1 correspondance with the
 connected components of the fibers of $q$. 
 
 Furthermore the natural map $f: S_{\infty} \to S_{N_0}$ has the following property: every point in $S_{N_0}$ has an open neighborhood $U$ such that
 each connected component of $f^{-1}(U)$ is finite over $U$. Since $S_{N_0}$ is Stein, it follows that $S_{\infty}$ is Stein, too  \cite{LeBa}. 
 \end{prv}

 \begin{rem}
  I don't know of a single connected quasi-K\"ahler manifold violating SSC(X). Note that SSC(X) implies (with $f=Id_X$) that if  $\pi_1(X)$ is infinite
 then $\pi_1(X)$ has  a finite dimensional complex representation with infinite image. 
 \end{rem}
\begin{ques}
 Is there a quasi-projective manifold $X$ with  an infinite fundamental group having no finite dimensional complex linear representation with infinite image? 
 Is there a quasi-projective manifold $X$ with  a simple infinite fundamental group?
\end{ques}

\noindent {Philippe Eyssidieux}\\
{Universit\'e  de Grenoble-Alpes. Institut Fourier.
100 rue des Maths, BP 74, 38402 Saint Martin d'H\`eres Cedex, France}\\
{philippe.eyssidieux@univ-grenoble-alpes.fr}
{http://www-fourier.ujf-grenoble.fr/$\sim$eyssi/}\\
\end{document}